\documentclass{amsart}

\usepackage{bm}
\usepackage{graphicx}
\usepackage{amsmath}
\usepackage{hyperref}
\usepackage{wrapfig}
\usepackage{amsfonts}
\usepackage{amsthm}

\numberwithin{equation}{section}
\newtheorem{thm}{Theorem}[section]

\newtheorem{remark}{Remark}
\newtheorem{definition}{Definition}

\newcommand{\be}{\begin{equation}}
\newcommand{\ee}{\end{equation}}
\newcommand{\la}{\label}

\begin{document}

\title[Effective distribution of codewords]{Effective distribution of codewords for Low Density Parity Check Cycle codes in the presence of disorder}

\author[Roshan Warman]{Roshan Warman}
\email{roshanwarman22@gmail.com}
\address{Academy at the Lakes}

\author[Iuliana Teodorescu]{Iuliana Teodorescu}
\email{iuliana@usf.edu}
\address{Department of Mathematics and Statistics, University of South Florida}

\author[Razvan Teodorescu]{Razvan Teodorescu}
\email{razvan@usf.edu}
\address{Department of Mathematics and Statistics, University of South Florida}



\begin{abstract}
We review the zeta-function representation of codewords allowed by a parity-check code based on a bipartite graph, and then investigate the effect of disorder on the effective distribution of codewords. The randomness (or disorder) is implemented by sampling the graph from an ensemble of random graphs, and computing the average zeta function of the ensemble. In the limit of arbitrarily large size for the vertex set of the graph, we find an exponential decay of the likelihood for nontrivial codewords corresponding to graph cycles. This result provides a quantitative estimate of the effect of randomization in cybersecurity applications. 
\end{abstract}

\maketitle

\section{Introduction}

\subsection{Motivation}
 In 1963, Robert Gallanger \cite{Gall} developed the Low Density Parity Check codes, which were powerful for their small size to efficacy ratio. In particular, an $(n, j, k)$ low density parity check matrix consists of $n$ block length, with $j$ ones in each column and $i$ ones in each row. Binary channels based on parity-check codes are described by a bipartite graph $(V, E)$ whose vertex set is the disjoint union of ``bits" $B$ and ``checks" $C$, $V = B \sqcup C$, such that edges $e \in E$ have one end in each subset: for $e = (v,w) \in E$ either $ (v, w) \in B \times C$ or $(v, w) \in C \times B$. A codeword allowed by this graphical structure is then any loop $\gamma = e_1 e_2 \ldots e_{2n}, \, \{e_k \} \subset E$. It corresponds to $n$ ``checks" (binary constraints) satisfied by a binary input of size $k \le n$. Generalizing this well-known construction, let $X(V,E)$ be an undirected graph and $\vec{X} (V, 2E)$ its associated directed graph, with directed adjacency matrix $M(\vec{X})$, of dimensions $2|E| \times 2|E|$ and entries $M_{ij} = 1$ if directed edge $e_i$ ``feeds" into directed edge $e_j$, and zero otherwise. The edge zeta function for the graph $X$ is defined by introducing a set of auxiliary (complex) variables for each edge, 
and defining $\zeta^{-1}_E(\{ u_k \},X) = \det [I - U\cdot M(\vec{X}])$, for $U = {\rm{ diag }}[u_1, u_2, \ldots u_{2|E|}], 
 \,\, u_{|E|+ k} = u_k.
$
The name is justified by the following product representation over the set of irreducible, non-backtracking, simple loops (or prime cycles) of the graph, $\mathcal{P}(X):$ $$\la{edge}
\zeta_E(\{ u_k \}, X) =  \prod_{\gamma \in \mathcal{P}(X)} \left [ 1-u_1 u_2 \ldots u_l \right ]^{-1}.
$$

The relation to Coding theory stems from the fact that $\zeta_E$ is a generating function for the pseudo-codewords of the graph (for cycle codes): $$
\zeta^{-1}_E(\{ u_k \}, X) = \sum_{n_1+ n_2 + \ldots \ge 0} C_{n_1 n_2 \ldots n_p} u_1^{n_1} u_2^{n_2} \ldots u_p^{n_p}, 
$$
where the coefficients $C_{n_1 n_2 \ldots n_p}$ are non-zero if and only if $(n_1, n_2, \ldots, n_p)$ 
is a pseudo-codeword of the graphical model $X$. 

\subsection{Goals and outline} The purpose of this study is to provide a characterization for the pseudo-codewords of Normal graphs randomly selected from an ensemble of parity check matrices. In particular, we present the ``average"value of these pseudo-codewords and with results from \cite{pseudo} present relevant applications to cybersecurity and the Shortest Vector Problem, effectively completing this endeavor. 

\section{Background}

\subsection{Pseudocodes of normal graphs}
Recall some basic notions of graph theory. A graph $G$ is defined by a triple $(V, E, \rho)$, where $V$ is the set of \textit{vertices}, $E$ is a set of \textit{edges}, which is disjoint from $V$, and $\rho$ is a function from $E \rightarrow P_2(V)$. A \textit{path} $(v_1, e_1, \dots, v_{n-1}, e_{n-1}, v_n)$ is a sequence of distinct vertices alternating with edges. A \textit{cycle} is a path where $v_n = v_1$. A cycle $\Gamma$ is \textit{backtrackless} if for no $i$, $e_i = e_{i+1}$. It is \textit{primitive} if there is no cycle $\Omega$, such that $\Gamma = \Omega^r$ (i.e., $\Gamma$ is obtained by traversing $\Omega$ $r \geq 2$ times), and it is \textit{equivalent} to $\Omega$ if there is some integer $p$ such that $e_i = e_{i \mod p}$.

Let $H$ be the parity check matrix of some binary linear code $\mathcal{C}$. Likewise, let $R(H)$ be the set of rows and $C(H)$ be the set of columns. We define $C_r(H) = \{c \in C(H) \mid H_{cr} = 1\}$, and similarly, $R_c(H) = \{r \in R(H) \mid H_{cr} = 1\}$.

\begin{definition}[Tanner graph]
Let $T(H)$ be the bipartite graph associated with the parity check matrix, $H$, with bit nodes, $X_1, \dots, X_{|C|}$, and check nodes, $p_1, \dots, p_{|R|}$ and bit node $i$ and check node $j$ connected if and only if $H_{ij} = 1$.
\end{definition}

\begin{remark}If $T(H)$ is $2$-regular for a given parity check matrix $H$, then we call the code a \textit{cycle-code} and its Tanner graph a \textit{Normal Graph}.  
\end{remark}
Our main results will be focused on these cycle codes, since defining codewords on acyclic codes is not only impractical, but also relatively trivial to decode.
In \cite{pseudo}, the authors make use of \textit{pseudo-codewords} in their development of decoding Low Density Parity Check codes with the zeta function. The cleverness of pseudo-codewords is to intentionally introduce a layer of redundancy in the graphical model which would make identifying unique elements easier as $M \rightarrow \infty$. In particular, the pseudo-codewords make use of \textit{unramified covers}, which is a graph homomorphism $\pi: X \rightarrow Y$. An \textit{$M$-cover} is an unramified cover such that for $x \in V(X),~ \pi^{-1}(X)$ contains $M$ vertices of $Y$. The advantage of decoding an $M$-fold becomes immediately clear as $M$ becomes larger.

\begin{definition}[Pseudo-codewords]
Let $C$ be a binary linear code with parity check matrix $H$ and Tanner graph, $T$. Further, let $\widetilde{T}$ be an $M$-fold cover of $T$. Then the pseudo-codeword of $\widetilde{c}$ is the vector $\bf{\omega}(\widetilde{c}) = (\omega_1(\widetilde{c}), \dots, \omega_n(\widetilde{c}))$, where $$ \omega_i(\widetilde{c}) = \dfrac{1}{M}\sum_{k \leq M} \widetilde{c}_{(i, k)}$$
and $\widetilde{c}$ is a codeword of the code $\widetilde{C}$ associated with $\widetilde{T}$.
\end{definition}

\begin{remark}
The pseudo-codewords exist in $\mathbf{Q}$, so the lifting of a codeword to its pseudo-codewords, while not unique, is surjective.
Also note that $\omega(\widetilde{c})\mod 2$ is a codeword of $C$.
\end{remark}

\subsection{Ihara Zeta function}
Our main tool in this endeavor will be the Ihara Zeta Function. The authors refer to \cite{Ihara} for a more complete review.  

\begin{definition}[Ihara Zeta Function]
The Ihara Zeta Function of $X$ is defined to be the power series $\zeta_X(u_1, \dots, u_n) \in \mathbb{Z}[u_1, \dots, u_n]$ given by
$$\zeta_{X}(u_1, \dots, u_n) = \prod_{\Gamma \in \mathcal{P}(X)}\left [ 1-u^{l(\gamma)}\right ]^{-1}, \quad u \in \mathbb{C}$$
\end{definition}

\begin{remark}
$\mathcal{P}(X)$ is the collection of equivalence classes of back trackless, tailless, and primitive cycles in $X$.
\end{remark}
Additionally, in 1989 Hashimoto-Bass proved a theorem which expresses the Ihara Zeta Function as a determinant:
\begin{thm}[Determinant-form Ihara Zeta Function]
Let $X$ be a graphical model, then $$\zeta(u, X)^{-1} = (1-u^2)^{r-1}\det(I - A + Qu^2),$$ where $r = |E| - |V| + 1.$
\end{thm}

The following theorem forms the basis of this project:

\begin{thm}[Koetter, et al. \cite{pseudo}]
Let $C$ be a cycle code defined by parity-check matrix $H$ having normal graph $T$, and let $\zeta_T(u_1, \dots, u_x)$ be its Ihara Zeta Function. Then the monomial $u_1^{p_1} \dots u_x^{p_x}$ has nonzero coefficients if and only if the corresponding exponent vector $(p_1, \dots p_x)$ is a pseudo-codeword for $C$. 
\end{thm}
\section{Main Results}

Consider now an ensemble of graphical models $\mathcal{G} = \{ X \}$, or (equivalently) an ensemble of directed adjacency matrices $\mathcal{M}$, of various sizes, described by a probability measure $d \mu$. We define the randomized zeta function of the ensemble as 
\be
\widetilde{\zeta}_{\mathcal{M}}(u) = \mathbb{E}_{M} [\det (I - u \cdot M)]^{-1} = \int_{\mathcal{M}} \det (I - u \cdot M)^{-1} d\mu(M),
\ee
%
where $u \in \mathbb{C}$ is a formal variable, and the definition is understood to be convergent in some convergence domain $|u| < R_{\mathcal{M}}$. The parameter $R_{\mathcal{M}} > 0$ will be referred to as convergence radius of the ensemble $\mathcal{M}$. 

\begin{remark}
In the random case, the averaged zeta function serves a formal power series whose coefficients provide quantitative estimates for the likelihood of loops corresponding to codewords, and therefore can be used to determine which random ensembles would be more vulnerable to attempts of finding a codeword at random. 
\end{remark}
To begin the analysis, note that the multivariate Gaussian provides an integral representation to the determinant of the resolvent of $M \in \mathbb{R}^{n \times n}$, by the identity
$$
\det(I - u M)^{-1}= 
\mathbb{E}_{Y, Z} \exp \left \{ \frac{u}{2} \text{Tr} [(YZ' +Z Y')M] \right \},
$$
where $Y, ~Z$ are two i.i.d. multivariate Gaussians in $\mathbb{R}^n$ with covariance matrix given by the identity. This leads to the representation 

\be 
\widetilde{\zeta}_{\mathcal{M}}(u) = \mathbb{E}_{Y, Z} \mathbb{E}_{M} \exp \left \{ \frac{u}{2} \text{Tr} [(YZ' +Z Y')M] \right \}
\ee
Throughout these formulas, the length of vectors $Y, Z$ is understood to correspond to the size of matrix $M$. As the random variables are independent $Y \perp Z \perp M \perp Y$, we obtain 

\be \la{formula}
\widetilde{\zeta}_{\mathcal{M}}(u)  = \mathbb{E}_{M} \sum_{k = 0}^{\infty} \frac{u^k}{2^k k!} \mathbb{E}_{Y, Z} \{ \text{Tr} [(YZ' +Z Y')M]\}^k
\ee

This representation presents a number of advantages when describing the distribution of codewords, as we indicate in the following subsection. 

\subsection{General properties of codeword distribution for randomized codes}

An immediate consequence of \eqref{formula} is that the power series contains only even powers of $u$, as all odd moments of the normalized Gaussian vanish. Therefore, the average of all ensembles allows only for graph cycles of even length, and the first nontrivial contribution is obtained by considering the coefficient of $u^4$ in \eqref{formula}: 
$$
\mathbb{E}_{Y, Z} (Tr [(YZ' +Z Y')M] )^4 =\sum_{i,j,k,l, p, q, r, s = 1}^{2n}  M_{ip}M_{jq} M_{kr} M_{ls} \mathbb{E}_{Y, Z} Y_iY_jY_kY_l 
Z_p Z_qZ_rZ_s
$$
up to an overall factor. Obviously, the result can be written as a linear combination over products of 2-cycles on one hand, and 4-cycles on the other, owing to the formulas
$$
\mathbb{E}_{Y} (Y_iY_jY_kY_l) = \delta_{ij}\delta_{kl} + \delta_{ik}\delta_{jl} + \delta_{il}\delta_{jk} + 3 \delta_{ij}\delta_{jk}\delta_{kl},
$$
where $\delta_{ij}$ is the Kronecker symbol. For example, the pairing of indices leading to the term 
$$
M_{ip}M_{qi} M_{kq} M_{pk} 
$$
corresponds to the 4-cycle $\gamma = e_ie_pe_ke_q$, and could be supported on two bit nodes and two check nodes (endpoints of edges $e_i, e_q$ and $e_q, e_k$ respectively). Analysis of the coefficients of higher-order terms in $u$ yields similar graphical interpretations. 

Evidently, the very distribution of random graphs (manifested by averaging over the ensemble of adjacency matrices $M$)
is the other contributing factor to the effective codeword distribution. In the absence of randomization, the relative weight of a codeword of length 4, for instance, could be easily estimated by truncating the formal power series \eqref{formula} to a polynomial, and using simple combinatorial identities to find the total maximal number of codeword configurations, and of 4-cycle codewords in particular. 

In order to describe the effect of randomization on the family of graphs $\{ X\}$, we need a few more theoretical tools, briefly reviewed in the next section.

\section{Large deviations principle extensions to graph theory} 
 
We recall the definition of a rate function for a random variable $Y$: it is is defined as
$$
I(x) := \max_{t > 0} \, \left [t x - \ln(m_Y(t)) \right ],
$$
where $m_Y(t)$ is the moment-generating function of $Y$, $m_Y(t) = \mathbb{E}(e^{tY}).$ The Large Deviations Principle  states that the probability of ``large deviations"
$$
P\left (\frac{1}{n}\sum_{k=1}^n X_k \ge x \right ) \sim e^{-n I(x)}, 
$$
where $X_1, X_2, \ldots, X_n$ is a sample of  i.i.d. r.v., and $I(x)$ is the rate function defined earlier, and $``\sim"$ means that the probability
is determined only up to an overall normalization factor.
For the case of interacting diffusions, assume that the sequence of random variables 
$Z_n, n = 1, 2, \ldots$ from the space $\Sigma$ have distributions 
$d P_n$ and moment-generating functions $m_n(t) = \mathbb{E}_{P_n}(e^{tZ_n})$.  
If the limit
\be \la{fenergy}
\Lambda(t) = \lim_{n \to \infty} n^{-1}\log m_n(t)
\ee
exists, is convex and bounded from below, and that 
\be \la{rate}
Q(x) = \sup_{t}[xt - \Lambda(t)]
\ee
is a well-defined {\emph{rate function}} (bounded from above,
lower semi-continuous and has compact level sets), then 
$\{ Z_n \}$ satisfies the Large Deviations Principle with rate  $Q$:
\be \la{one}
\begin{array}{c}
\lim \sup  \log P(Z_n \in C ) \le -n Q(C), \\ 
\lim \inf \log P(Z_n \in O ) \ge -n Q(O), \,\!\, 
\end{array}
\ee
where $n \to \infty$, the sets $C, O$ are closed and open, respectively, 
and $Q(S)$ is by definition 
\be
Q(S) \equiv \inf_{x \in S} Q(x), \,\, (\forall) \,\, S \subset \Sigma.
\ee
This is known as the G\"artner-Ellis theorem \cite{G-E, ge}. 
For many situations, (\ref{one}) imply that for any set $S$, 
\be \la{ge}
n^{-1} \log P(Z_n \in S) \to - Q(S).
\ee
 
 In \cite{Var}, the authors extended this principle to the estimation of the number of loops of given length in a random graph from the Erd\"os-R\'enyi class $G(N, p)$ (where $N$ stands for the number of vertices, and $p$ stands for the probability to establish an edge between two vertices), in the limit $N \to \infty$. Specifically, they computed a rate function $\phi_p(t)$ for the random variable $T_{N, p}$ (the number of cycles of length 3 in any realization of the graph $X \in G(N, p)$), and found (Theorem 4.1, \cite{Var}) the scaling 
 $$
 \mathbb{P}(T_{N,p} \ge N^3 t) \sim e^{-N^2 \phi_p(t)}, \, t > 0, \, p \in (0, 1),
 $$ 
 where $\phi_p(t) = O(1)$ with respect to asymptotic growth in $N$. By a similar argument (to be presented in a forthcoming publication), it is possible to obtain the scaling 
 $$
 \mathbb{P}(L_{N,p} \ge N^4 t) \sim e^{-N^2 \varphi_p(t)}, \, t > 0, \, p \in (0, 1),
 $$ 
with a different rate function $\varphi_p(t)$, and where $L_{N, p}$ stands for the number of 4-cycles. This indicates a very fast exponential suppression of 4-cycles (conjecturally, of all short loops, i.e. loops of length $O(N^0)$), and, correspondingly, of codewords in the randomized code ensemble. This fact has obvious important applications in cybersecurity. A quantitative analysis of this type of exponential suppression will be provided separately. 

\subsection{Randomization with protected subspaces} 

A variant of the complete randomization considered here, most relevant for applications, consists of embedding an invariant subgraph with short loops into an ensemble of random graphs. This can be realized formally by generating all conjugacy classes of a given adjacency matrix, and then sampling from the entire group with uniform probability. The difficulty is, of course, that a generic element from a  group of conjugacy classes does not have a meaningful interpretation as adjacency matrix of some graph. However, the problem of determining (by random sampling) whether a group of matrices has been generated by conjugacy of a set of adjacency matrices (the ``hacker's problem") is interesting and nontrivial in its own right; it bears more than just a superficial similarity to both the Hidden Subgroup Problem in cryptography, and to the invariant subspace problem in function theory, and it may be regarded as a probabilistic variant to both. We will address this problem and its implications in a future publication.

 \section{Appendix}
 \url{https://github.com/Roshanwarman/LDPC-Zeta-function-decoding}
 
 \vspace{2mm}

\end{document}